\documentclass[12pt]{elsarticle}

\usepackage{amssymb}

\journal{arXiv.org}

\newcommand{\grad}{\mathop{\rm grad}\nolimits}
\renewcommand{\div}{\mathop{\rm div}\nolimits}

\begin{document}

\begin{frontmatter}

\title{VAGO method for the solution of elliptic second-order 
boundary value problems}

\author{Nikolay P. Vabishchevich}
\ead{nvab@ibrae.ac.ru}
\address{Nuclear Safety Institute,
52, B. Tulskaya, 115191 Moscow, Russia}

\author{Petr N. Vabishchevich}
\ead{vabishchevich@gmail.com}
\address{Keldysh Institute of Applied Mathematics, 
4-A Miusskaya Square, 125047 Moscow, Russia}

\begin{abstract}

Mathematical physics problems are often formulated using differential oprators of vector analysis -
invariant operators of first order, namely, divergence, gradient and rotor operators.
In approximate solution of such problems it is natural to employ similar operator formulations
for grid problems, too. The VAGO (Vector Analysis Grid Operators) method is based on such a
methodology.
In this paper the vector analysis difference operators 
are constructed using the Delaunay 
triangulation and the Voronoi diagrams. 
Further the VAGO method is used to solve approximately boundary value problems
for the general elliptic equation of second order.
In the convection-diffusion-reaction equation the diffusion coefficient is
a symmetric tensor of second order.

\end{abstract}


\end{frontmatter}

\section{Introduction}
\label{sec:into}

In the theory of difference schemes \cite{Samarskii_theory} the problem
of approximation of a differential problem via grid one received primary emphasis.
A universal method of balance (integro-interpolation method)
is used as the basic one to construct discrete analogs.
It was proposed by Samarskii \cite{balansametod} and currently it is known as
the control volume method \cite{FVM,Finite_volume_methods,Generalized_difference_methods}.
A difference scheme is constructed using the balance method via integration
of the initial equation over a control volume --- a part of the computational
domain adjacent to the specified grid point. 

In solving mathematical physics problems on general unstructured grids
there is often used the Delaunay triangulation which displays 
some optimal features \cite{George_Louis_Borouchaki}. 
For the Delaunay triangulation it is natural to consider as control volumes
a Voronoi polygon --- a set of points lying closer to the considered
grid point in compare with all others.
Consistent approximations of the vector analysis operators on the basis of 
the Delaunay triangulation were constructed in \cite{DirichleFavorskii} for
some problems of mathematical physics. A systematic investigation of approxinations
for gradient and divergence operators on dual Voronoi-Delaunay partitionings
has been conducted in \cite{Nicolaides_Wang,Nicolaides_Wu} for model 2D and 3D 
problems.

The present work continues our studies on Vector Analysis Grid Operators method.
In paper \cite{vago} there are discussed problems of constructing
consistent approximations for basic operators of vector analysis ---
gradient, divergence amd rotor operators. Special representations were derived
for truncation errors that allowed to obtain the corresponding estimates
for the convergence rate of the approximate solution to exact one.
Possibilities of the VAGO method have been demonstrated on some scalar 
and vector problems. Difference schemes for steady-state convection-diffusion
problems were constructed on the basis of the dual Voronoi-Delaunay partitioning
in \cite{Mishev,Vabishchevich_Samarskij2000}. In these problems 
(see, e.g., \cite{Lazarov_Vassilevski}) the emphasis is on properties
of discrete operators for the convection and diffusion transport.
The convective terms are written \cite{samvabkon,Vabishchevich_Samarskij2000}
in the divergent (conservative), non-divergent (characterictic) and 
symmetric forms. In \cite{vago} there were designed VAGO
approximations as well as investigated truncation errors 
for the Dirichlet problem for stationary and usteady
convection-diffusion problems for isotropic media. 

Paricular emphasis should be placed on approximation of
diffusion transport operator for an anisotropic medium
where the diffusion coefficient is a symmetric tensor of second order.
In this case we should control at the discrete level not only
the seld-adjoint property of the grid diffusion operator but the 
monotonicity and conservation properties for the corresponding grid 
approximatioms, too \cite{Samarskii_theory}.
The present-day view on approximation of problems with mixed derivatives
on rectangular grids is given in \cite{Matus_Rybak}. 
Properties of standard finite-element approximations on triangular or 
tetrahedral meshes are discussed in \cite{Putti_Cordes} in detail.
To obtain the monotonicity property, a modification of standard 
linear Galerkin finite element discretizations of the Laplace operator
was performed. Much attention is given to constructing control volume
schemes for the approximate solution of convection-diffusion
boundary value problems with an anisotropic diffusion coefficient
(see \cite{Droniou_Eymard} and references therein).
In particular (see also \cite{Eymard_Gallouet_Herbin}),
the control volume schemes are developed with an anisotropic diffusion
for approximating a local discrete gradient.
The results of work \cite{Nicolaides_Wang} are generalized 
in \cite{Hu_Nicolaides} to the 2D div-curl system in anisotropic
media.

The outline of this paper is the following.
In Section 2, the boundary-value problem for the general
elliptic equation of second order is formulated in an invariant form
as the problem for a convection-diffusion equation based on
using vector analysis operators --- divergence and gradient.
The diffusion coefficient is a symmetric positive definite tensor
of second order. On dual Voronoi-Delaunay partitionings (Section 3)
there are introduced the corresponding spaces of grid functions.
Using the VAGO method the difference operators of divergence and gradient
are constracted and investigated in Section 4. Difference schemes for
the considered convection-diffusion-reaction problem are studied in 
Section 5. Approximation of the diffusion transport operator with an
anisotropic diffusion coefficient has reseived much consideration.

\section{The boundary value problem}

The boundary value problem is considered for an elliptic equation 
of second order. In $n$-dimensional bounded  domain $\Omega$
with smooth enough boundary $\partial \Omega$ function 
$u(\mathbf{x}), \ \mathbf{x} = (x_1, ..., x_n)$ satisfies 
the following equation
\[
- \sum_{{\alpha},\,{\beta}=1}^n \frac{\partial }{\partial
x_{\alpha}} \left( d_{{\alpha}{\beta}}(\mathbf{x}) \frac{\partial u}{\partial
x_{\beta}} \right) + 
\sum_{{\alpha}=1}^n \frac{\partial }{\partial x_{\alpha}} (a_{\alpha} \ u) +
\sum_{{\alpha}=1}^n b_{\alpha} \frac{\partial u}{\partial x_{\alpha}} +
\]
\begin{equation}\label{1}
c(\mathbf{x}) u(\mathbf{x}) = f(\mathbf{x})  , \quad
\mathbf{x} \in \Omega .
\end{equation}
The coefficients at the higher order derivatives meet the ellipticity
condition:
\[
    d_{{\alpha}{\beta}} = d_{{\beta}{\alpha}} , 
    \quad \sum_{{\alpha},\,{\beta}=1}^n
    d_{{\alpha}{\beta}}(\mathbf{x}) \, y_{\alpha} \, y_{\beta} \geq
    \delta  \, \sum_{{\alpha}=1}^n y_{\alpha}^2 , 
\]
\begin{equation}\label{2}
	\quad \delta >0, \quad
    \forall \, \mathbf{y} \in \mathbb{R}^n, \quad \forall \,
    \mathbf{x} \in \overline \Omega .
\end{equation}
For other coefficients of equation (\ref{1})  we have
\begin{equation}\label{3}
  \frac{1}{2} \sum_{{\alpha}=1}^n 
  \frac{\partial (a_{\alpha} - b_{\alpha})}{\partial x_{\alpha}} + 
  c(\mathbf{x}) \ge 0 .
\end{equation}
Assume that the right hand side of the equation can be represented 
as follows
\begin{equation}\label{4}
  f(\mathbf{x}) = f_0(\mathbf{x}) +
  \sum_{{\alpha}=1}^n \frac{\partial f_{\alpha}}{\partial
  x_{\alpha}} .
\end{equation}
The simplest boundary value problem is considered for (\ref{1}) with
the following boundary conditions
\begin{equation}\label{5}
u(\mathbf{x}) = 0, \quad \mathbf{x} \in \partial \Omega .
\end{equation}

In Hilbert space $L_2(\Omega)$ the scalar product and norm 
are introduced like this
\[
  (u, v) = \int_{\Omega} u(\mathbf{x}) v(\mathbf{x}) d \mathbf{x},
  \quad \|u\| = (u, u)^{1/2} .
\]
In $\mathcal{H}(\grad)$ we introduce
\[
  (u, v)_{\mathcal{H}(\grad)} = (u,v) + 
  \sum_{{\alpha}=1}^n \left ( \frac{\partial u}{\partial x_{\alpha }},
  \frac{\partial u}{\partial x_{\alpha }} \right ) ,
  \quad \|u\|_{\mathcal{H}(\grad)} = \left ( (u, u)_{\mathcal{H}(\grad)}\right )^{1/2} .
\]
For the solution of problem (\ref{1}) - (\ref{5}) the following apriori 
estimate takes place (see, e.g., \cite{Gilbarg})
\begin{equation}\label{6}
    \|u\|_{\mathcal{H}(\grad)} \leq M \left (\sum_{{\alpha}=0}^n
    \|f_{\alpha}\|^2 \right )^{1/2}.
\end{equation}

Unstructured computational grids are often used to solve approximately
boundary problems of mathematical physics. In some cases it is not
suitable to use the coordinate-wise form of the considering
equations and boundary conditions in one or another coordinate system.
It seems more reasonable to employ the invarinat formilation of
the problem. In this case the considered problem is represented
via operators of vector and tensor analysis.

Equation (\ref{1}) is the basic one for problems of
continuum mechanics. Vectors
\[
  \mathbf{a} = \{a_{\alpha} \},
  \quad \mathbf{b} = \{b_{\alpha} \},
  \quad \alpha =1,... , n,
\]
are connected with the convective transport.
Symmetric tensor of second rank
\[
  \mathbf{D} = \{d_{{\alpha}{\beta}} \},
  \quad \alpha =1,... , n,
  \quad \beta  =1,... , n 
\]
in (\ref{1}) corresponds to the diffusion transport.
Coefficient $c(\mathbf{x})$ can be treated as a reaction.
In such an interpretation equation (\ref{1}) is nothing
buth the equation of convection-diffusion-reaction
\begin{equation}\label{7}
  - \div (\mathbf{D}\grad \,u) + \div (\mathbf{a} u) +
  \mathbf{b} \grad \,u 
  + c(\mathbf{x}) u(\mathbf{x}) = f(\mathbf{x})  , \quad
  \mathbf{x} \in \Omega .
\end{equation}
Boundary value problem (\ref{1}) - (\ref{5}) in form (\ref{5}), (\ref{7})
does not employ explicitly Cartesian coordinates
$x_1, ..., x_n$.

\section{Delaunay triangulation and Voronoi diagrams}

Assume that the computational domain is a convex polyhedron $\Omega, \ n =2,3$ with 
boundary $\partial\Omega$. In domain 
$\overline{\Omega} = \Omega \cup \partial\Omega$ we consider the 
grid $\overline{\omega}$, which consists of nodes 
$\mathbf{x}^{D}_{i}, \ i = 1,2,\ldots,M_D$, and the angles of the polyhedron 
$\Omega$ are nodes. Let $\omega$ be a set of inner nodes and 
$\partial \omega$ is a set of boundary nodes, i.e., 
$\omega=\overline{\omega} \cap \Omega$, 
$\partial \omega = \overline{\omega} \cap \partial \Omega$.

Each node $\mathbf{x}^{D}_i, \ i = 1,2,\ldots,M_D,$ connect a
certain part of the computational domain, namely, the Voronoi polyhedron or 
its part belonging to $\Omega$. The Voronoi polyhedron (polygon) for a separate 
node is a set of points lying closer to this node than to all the other ones:
\[
  V_i = \{ \, \mathbf{x} \ | \ \mathbf{x} \in \Omega,
   \ | \mathbf{x} - \mathbf{x}^D_i | < | \mathbf{x} - \mathbf{x}^D_j |, \
  j = 1,2,\ldots, M_D \, \},
  \quad i = 1,2,\ldots, M_D ,
\]
where $|\cdot|$ is the Euclidean distance. Each vertex 
$\mathbf{x}^{V}_k, \ k = 1,2,\ldots, M_V$ of the Voronoi polyhedron is 
associated with the tetrahedron constructed by the appropriate nodes 
contacting the Voronoi polyhedrons. We will assume that all vertices 
of the Voronoi polyhedrons lie either inside the computational domain $\Omega$ or 
on its boundary $\partial \Omega$. These tetrahedrons determine the Delaunay 
triangulation --- a dual triangulation to the Voronoi diagram. The $D$-grid in the domain 
$\Omega$ is determined by the set of nodes (vertices of tetrahedrons of the 
Delaunay triangulation) $\mathbf{x}^{D}_i, \ i = 1,2,\ldots,M_D$, the $V$-grid is defined 
by the set of nodes (vertices of polyhedron of the Voronoi diagram) 
$\mathbf{x}^{V}_k, \ k = 1,2,\ldots,M_V$.

We mark a separate tetrahedron $D_k$ of the Delaunay triangulation.
This tetrahedron is identified by the number $k$ of 
the Voronoi polyhedronm vertex, $k = 1,2, \ldots,M_V$. The 
tetrahedrons $D_k, k = 1,2, \ldots,M_V$ cover the entire computational domain, so
\[
  \overline{\Omega} = 
  \mathop{\cup}^{~M_V}_{k=1} 
  \overline{D}_{k},
  \ \overline{D}_{k} = D_{k} \cup \partial D_{k},
  \ D_{k} \cap D_{m} = \emptyset,
  \ k \neq m,
  \ k, \ m = 1,2,\ldots,M_V.
\]
For common planes of the tetrahedron we use the notations
\[
  \partial D_{km} = \partial D_{k} \cap \partial D_{m} ,
  \quad k \neq m,
  \quad k, \ m = 1,2,\ldots,M_V.
\]
The boundary of the computational domain $\partial \Omega$ consists of the 
planes of Delaunay tetrahedrons. Let
\[
  \partial D_0 = \partial \Omega,
  \quad 
  \partial D_{k0} = \partial D_{k} \cap \partial D_{0} ,
  \quad k = 1,2,\ldots,M_V.
\]

We associate the Voronoi polyhedron $V_i, \ i = 1,2,\ldots, M_D$ with the node of the main grid $i$.
Thus, we have
\[
  \overline{\Omega} = 
  \mathop{\cup}^{~M_D}_{i=1} 
  \overline{V}_{i},
  \quad \overline{V}_{i} = V_{i} \cup \partial V_{i},
  \quad V_{i} \cap V_{j} = \emptyset,
  \quad i \neq j,
  \quad i, \ j = 1,2,\ldots,M_D
\]
and 
\[
  \partial V_{ij} = \partial V_{i} \cap \partial V_{j} ,
  \quad i \neq j,
  \quad i, \ j = 1,2,\ldots,M_D.
\]

For each tetrahedron $D_k$,  $k = 1,2, \ldots,M_V$, we define the set of 
neighbors $\mathcal{W}^{D}(k)$, having common planes with $D_k$, i.e., 
\[
  \mathcal{W}^{D}(k) = \{ m \ | \ \partial D_{k} \cap \partial D_{m} \neq \emptyset,
  \ m =0,1,\ldots,M_V \},
  \quad k = 1,2,\ldots,M_V.
\]
In this case, $m=0$ means that the tetrahedron $D_k$ contacts the boundary.
We define also the neighbors for each Voronoi polyhedron
$V_i, \ i = 1,2,\ldots, M_D$:
\[
  \mathcal{W}^{V}(i) = \{ j \ | \ \partial V_{i} \cap \partial V_{j} \neq \emptyset,
  \ j =1,2,\ldots,M_D\},
  \quad i = 1,2,\ldots,M_D.
\]

We assume that the introduced Delaunay triangulation and the Voronoi diagram are regular
\cite{PGCiarlet_1978a}. For the notations
\[
  h^{D}_k = \mathrm{diam} (D_k) \textrm{ --- diameter }D_k,
\]
\[
  \varrho^{D}_k = \mathrm{sup} \, \{\mathrm{diam}(S) \ | \ 
  S \textrm{ --- sphere in } D_k \} ,
  \quad k = 1,2,\ldots,M_V
\]
the regularity condition of the Delaunay triangulation is
\[
  \frac{h^{D}_k}{\varrho^{D}_k} \le \sigma  > 0, \quad k = 1,2,\ldots,M_V.
\]
Likewise, for the Voronoi diagram we have
\[
  h^{V}_i = \mathrm{diam} (V_i) \textrm{ --- diameter } V_i,
\]
\[
  \varrho^{V}_i = \mathrm{sup} \, 
  \{\mathrm{diam}(S) \ | \ S \textrm{ --- sphere in } V_i \} ,
\]
\[
  \frac{h^{V}_i}{\varrho^{V}_i} \le \sigma  > 0, \quad i = 1,2,\ldots,M_D,
\]
and
\[
  h = \max_{i,k} \{h_i^{V}, h^{D}_k \},
\]
\[
  \mathrm{meas} (D_k) = \int_{D_k} d \mathbf{x},
  \quad k = 1,2,\ldots,M_V ,
\]
\[
  \mathrm{meas} (V_i) = \int_{V_i} d \mathbf{x},
  \quad i = 1,2,\ldots,M_D .
\]

We will approximate the scalar functions of the continuous argument by the scalar grid functions 
that are defined in the nodes of the $D$-grid or in the nodes of the $V$-grid. We denote 
by $H_{D}$ the set of grid functions defined on the $D$-grid
\[
  H_{D} = \{ \ y(\mathbf{x}) \ | \ y(\mathbf{x}) = y(\mathbf{x}^{D}_i)
  = y^{D}_i, \ i = 1,2,\ldots,M_D \, \} .
\]
For the functions $y(\mathbf{x}) \in  H_{D}$, vanishing on the boundary $\partial \omega$,
we define 
\[
  H^{0}_{D} = \{ \ y(\mathbf{x}) \ | \ y(\mathbf{x}) \in H_{D}, \
  y(\mathbf{x}) = 0, \ \mathbf{x} \in \partial \omega  \, \} .
\]
We consider the scalar product and the norm for the scalar grid functions from $H_{D}$ by 
\[
  (y, v)_{D} =
  \sum_{i=1}^{M_D} 
  y_{i}^{D} \, v_{i}^{D} \ 
  \mathrm{meas} (V_{i}),
  \quad \|y\|_D = (y, y)_{D}^{1/2}.
\]
This scalar product and the norm are grid analogs of the scalar product and the
$L_2(\Omega)$-norm for the scalar functions of the continuous argument.

To determine the vector field in the control volume, it is natural to use 
the components of the sought function normal to the corresponding planes 
of the control volume. Choosing the initial and the final node, we connect 
with each tetrahedron edge $D_k, \ k = 1,2,\ldots,M_V$ or polyhedron edge 
$V_i, \ i = 1,2,\ldots,M_D$, a vector --- the directed edge. For Delaunay 
triangulation the normals to the planes are the directed edges of Voronoi diagram 
and vice versa. For the approximation of the vector 
functions thereby we can use projections of the vectors on the directed edges. We will 
further use exactly this variant with the description of the vector field in the 
control volume by means of vector projections on the edges of the control volume.

We will orient the Delaunay triangulation edges by the unit vector
\[
  \mathbf{e}^{D}_{ij} = \mathbf{e}^{D}_{ji},
  \quad i = 1,2,\ldots,M_D, \quad j \in \mathcal{W}^{V}(i),
\]
directed from the node with a smaller number to the node of a larger number. 
The vector function $\mathbf{y}(\mathbf{x})$ on the Delaunay triangulation is 
defined by the components
\[
  y_{ij}^{D} = \mathbf{y} \, \mathbf{e}^{D}_{ij},
  \quad i = 1,2,\ldots,M_D , \quad j \in \mathcal{W}^{V}(i),
\]
that are given in the middle of the edges
\[
  \mathbf{x}^{D}_{ij} = \frac{1}{2}(\mathbf{x}^{D}_i + \mathbf{x}^{D}_j).
\]
For the Delaunay triangulation used and the Voronoi diagram we define the length of 
the edges in the following way:
\[
  l^{D}_{ij} = | \mathbf{x}^{D}_{i} - \mathbf{x}^{D}_{j} |,
  \quad i = 1,2,\ldots,M_D , \quad j \in \mathcal{W}^{V}(i). 
\]

We denote by $\mathbf{H}_D$ the set of grid vector functions determined by 
the components $y_{ij}^{D}, \ i = 1,2,\ldots,M_D, \quad j \in \mathcal{W}^{V}(i)$ 
that are given in the middle of the edges. In a similar way we denote 
by $\mathbf{H}_V$ the set of grid vector functions defined by the 
components $y_{km}^{V}, \ k = 1,2,\ldots,M_V, \quad m \in \mathcal{W}^{D}(k)$. 
If the tangential components of the grid vector functions $\mathbf{y} \in \mathbf{H}_D$ 
vanish on the boundary, we define 
\[
  \mathbf{H}_D^0 = \{ \, \mathbf{y} \ | \ \mathbf{y} \in \mathbf{H}_D, 
  \quad \mathbf{y} (\mathbf{x}) \, \mathbf{e}^{D}_{ij}  = 0,
  \quad \mathbf{x} = \mathbf{x}^{D}_{ij} \in \partial \omega,
\]
\[
  \quad i = 1,2,\ldots,M_D, \quad j \in \mathcal{W}^{V}(i) \, \} .
\]
Consider the scalar product and norm 
\[
  (\mathbf{y}, \mathbf{v})_{D} =
  \sum_{i=1}^{M_D} 
  \sum_{j \in \mathcal{W}^{V}(i)} \
  y_{ij}^{D} \, v_{ij}^{D} \ 
  |\mathbf{x}^{D}_i - \mathbf{x}^{D}_{ij} | \ \mathrm{meas} (\partial V_{ij}),
  \quad 
  \|\mathbf{y}\|_D = (\mathbf{y}, \mathbf{y})_{D}^{1/2} .
\]
in $\mathbf{H}_D$.

\section{Gradient and divergence approximations}

Let us present the main results (see \cite{vago} for deatails)
connected with approximation of gradient and divergence operators
on the basis of the Delaunay triangulation and Voronoi diagram.

The set of grid functions $H_D$ can be the domain of definition of the 
grid gradient operator. We denote the grid gradient operator
by $\mathrm{grad}_D$. Taking into account
the chosen edge orientation, at the points $\mathbf{x}^{D}_{ij}$ we set
\begin{equation}\label{8}
  (\mathrm{grad}_D y)^D_{ij} = \eta(i,j) \ \frac{y_j^D - y_i^D}{l_{ij}},
  \quad i = 1,2,\ldots,M_D , \quad j \in \mathcal{W}^{V}(i),
\end{equation}
i.e., the range of values of the operator $\mathrm{grad}_D: \, H_D \rightarrow \mathbf{H}_D$
is the set of vector grid functions $\mathbf{H}_{D}$. In (\ref{8}), we use 
the following notation:
\[
  \eta(i,j) = \left \{
\begin{array}{ll}
  1,  & j > i,  \\
  -1,  & j < i .  \\
\end{array}
\right .
\]
For the truncation error of the grid operator $\mathrm{grad}_D$
we have
\begin{equation}\label{9}
  (\mathrm{grad}_D u) ({\mathbf{x}}) =
  (\mathrm{grad} \, u) ({\mathbf{x}})  + 
  \mathbf{g}(\mathbf{x}),
  \quad \mathbf{g} = \mathcal{O}(h^2), 
  \quad \mathbf{x} = \mathbf{x}_{ij}^{D},
\end{equation}
\[
  \quad i = 1,2,\ldots,M_D , \quad j \in \mathcal{W}^{V}(i)
\]
provided that $u(\mathbf{x})$ is a sufficiently smooth functions.

For the Voronoi polyhedron 
this equality is written in the following form:
\begin{equation}\label{10}
  \int\limits_{V_i} \mathrm{div} \, \mathbf{u} \ d \mathbf{x} =
  \sum_{j \, \in \mathcal{W}^{V}(i)} \ \int\limits_{\partial V_{ij}}
  (\mathbf{u} \cdot \mathbf{n}^{V}_{ij}) \ d \mathbf{x}, 
\end{equation}
where $\mathbf{n}^{V}_{ij}$ is the normal to the edge $\partial V_{ij}$ outside 
with respect to  $V_i$. To construct the grid operator $\mathrm{div}_D: \, 
\mathbf{H}_D \rightarrow H_D$ we use the elementary formulas of integration 
for the left- and right-hand sides (\ref{10}).
This leads to the grid analogs of the operator 
$\mathrm{div}$ in the form of
\begin{equation}\label{11}
  (\mathrm{div}_{D} \, \mathbf{y} )_{i}^D =
  \frac{1}{\mathrm{meas}(V_{i})}
  \sum_{j \, \in \mathcal{W}^{V}(i)} \, 
  (\mathbf{n}^{V}_{ij} \cdot \mathbf{e}_{ij}^{D}) \,
  y_{ij}^{D} \,
  \mathrm{meas}(\partial V_{ij}) ,
  \quad i = 1,2,\ldots,M_D .
\end{equation}
We have
\[
  (\mathrm{div}_D \mathbf{u}) ({\mathbf{x}}) =
  (\mathrm{div} \, \mathbf{u}) ({\mathbf{x}})  + 
  g(\mathbf{x}) +
  (\mathrm{div}_D \mathbf{q}) ({\mathbf{x}}),
\]
\begin{equation}\label{12}
  \quad g = \mathcal{O}(h), 
  \quad \mathbf{q} = \mathcal{O}(h), 
  \quad \mathbf{x} = \mathbf{x}_{i}^{D},
  \quad i = 1,2,\ldots,M_D ,
\end{equation}
for the truncation error. 
The truncation error for the grid divergence operators
$\mathrm{div}_D$ equals to $\mathcal{O}(1)$. 
However, there does exist a special divergence expression of the 
truncation error saving the situation in the case of approximation of 
problems of mathematical physics.

Let us highlight the adjoint property of the considered grid operators of
gradient and divergence. Taking into account the above notations we have
\[
  (v, \mathrm{div}_{D} \, \mathbf{y})_D +
  (\mathrm{grad}_{D} \, v, \mathbf{y})_D =
\]
\begin{equation}\label{13}
  \sum_{i \in \partial \omega } \, 
  \sum_{j \, \in \mathcal{W}^{V}(i), \ j \in \partial \omega } \, 
  (\mathbf{n}^{V}_{ij} \cdot \mathbf{e}_{ij}^{D}) \,
  \frac{v^D_j + v^D_i}{2} \,
  y_{ij}^{D} \,
  \mathrm{meas}(\partial V_{ij}) .
\end{equation}
From (\ref{13}) it follows that
\begin{equation}\label{14}
  \mathrm{div}^*_D = - \, \mathrm{grad}_D 
\end{equation}
on the set of grid functions $v \in H_D^0$ and $\mathbf{y} \in \mathbf{H}_D$.

\section{Approximation of the convection-diffusion-reaction problem}

In domain $\Omega$ we consider the regular $D$-grid and $V$-grid. 
Taking into account  boundary condition (\ref{5}), we will find the approximate solution of 
problem (\ref{5}), (\ref{7}) as grid function $y(\mathbf{x}) \in H_D^0$. 

A particular discussion should be given to approximation of convective,
 and moreover, diffisive terms in equation (\ref{7}).

We consider approximations of the convective transport operator both
in the divergent and in non-divergent (characteristic) form:
\[
  \mathcal{C}_1(\mathbf{v}) u = \mathbf{v} \grad u,
  \quad 
  \mathcal{C}_2(\mathbf{v}) u = \div (\mathbf{v} u) 
\]
at specified vector field $\mathbf{v}(\mathbf{x}), \ \mathbf{x} \in \Omega$.

The tangential components are approximated at the edge midpoints of 
the Delaunay triangulation as follows
\[
  a_{ij}^{D} =
  \frac{1}{ \mathrm{meas}(\partial V_{ij}) } \,
  \int\limits_{\partial V_{ij}} 
  \mathbf{a}(\mathbf{x}) \mathbf{e}_{ij}^{D} \, d \mathbf{x},
  \quad i = 1,2,\ldots,M_D, \quad j \in \mathcal{W}^{V}(i) .
\]
For scalar functions it is natural to introduce
\[
  y(\mathbf{x}^D_{ij}) =
  \frac{1}{2} \, (y(\mathbf{x}^D_{i}) + y(\mathbf{x}^D_{j})) .
\]
The direct approximation of the convective transport operator in the divergent form
results in
\begin{equation}\label{15}
  \mathcal{C}_2^D (\mathbf{v}^D) y = \div_D (\mathbf{v}^D y) =
  \frac{1}{2 \mathrm{meas}(V_{i})}
  \sum_{j \, \in \mathcal{W}^{V}(i)} \, 
  (\mathbf{n}^{V}_{ij} \cdot \mathbf{e}_{ij}^{D}) \, v_{ij}^{D} \,
  (y_{j}^{D} + y_{i}^{D}) \,
  \mathrm{meas}(\partial V_{ij}) .
\end{equation}
As for the convective transport operator in the non-divergent form,
its approximation is performed \cite{Vabishchevich_Samarskij2000}
using the grid analog of the following equality 
\[
  \mathcal{C}_1(\mathbf{v}) u = \mathcal{C}_2(\mathbf{v}) u - u \div \mathbf{v} .
\]
At
\[
  \mathcal{C}_1^D (\mathbf{v}^D) y = \mathcal{C}_2^D (\mathbf{v}^D) y - 
  y \div_D \mathbf{v}^D 
\]
we have
\begin{equation}\label{16}
  (\mathcal{C}_1^D (\mathbf{v}^D)\, \mathbf{y} )_{i}^D =
  \frac{1}{2 \mathrm{meas}(V_{i})}
  \sum_{j \, \in \mathcal{W}^{V}(i)} \, 
  (\mathbf{n}^{V}_{ij} \cdot \mathbf{e}_{ij}^{D}) \, v_{ij}^{D} \,
  (y_{j}^{D} - y_{i}^{D}) \,
  \mathrm{meas}(\partial V_{ij}) .
\end{equation}
For the trancation error similarly to (\ref{12}) we have
\[
  \mathcal{C}_{\alpha}^D (\mathbf{v}^D) u = 
  \mathcal{C}_{\alpha} (\mathbf{v}) u +  
  g(\mathbf{x}) +
  (\mathrm{div}_D \mathbf{q}) ({\mathbf{x}}),
  \quad \alpha =1,2, 
\]
\begin{equation}\label{17}
  \quad g = \mathcal{O}(h), 
  \quad \mathbf{q} = \mathcal{O}(h), 
  \quad \mathbf{x} = \mathbf{x}_{i}^{D},
  \quad i = 1,2,\ldots,M_D .
\end{equation}

Let us construct a grid analog of the diffusion transport operator
\[
  \mathcal{D} u = 
  \div (\mathbf{D}\grad \,u) , \quad
  \mathbf{x} \in \Omega .
\]

\begin{figure}
  \begin{center}
    \includegraphics[scale=1]{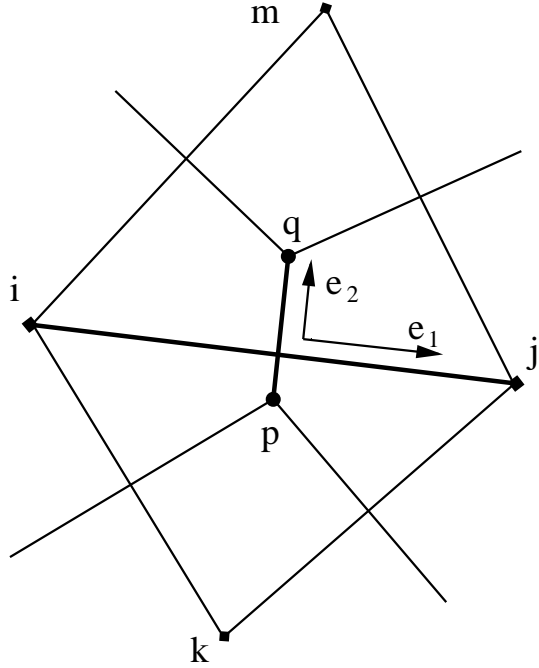} 
	\caption{Two-dimensional case}
    \label{f1}
  \end{center}
\end{figure}

We put formally
\begin{equation}\label{18}
  \mathcal{D}^D y = 
  \div_D (\mathbf{D}^D\grad_D y) .
\end{equation}
Above (see (\ref{8})) we presented the grid approximation for the projection
of the gradient operator onto edges of the Delaunay triangulation. 
It is enough for approximation of the diffusion operator in isotropic media
(with a scalar diffusion coefficient). For general elliptic equations
with mixed derivatives (\ref{1})  (with a tensor diffusion coefficient)
it is necessary to approximate the gradient not only along the edges of triangulation.
To investigate this problem, let us start from the 2D case ($n=2$, Fig.\ref{f1}).

\begin{figure}
  \begin{center}
    \includegraphics[scale=1]{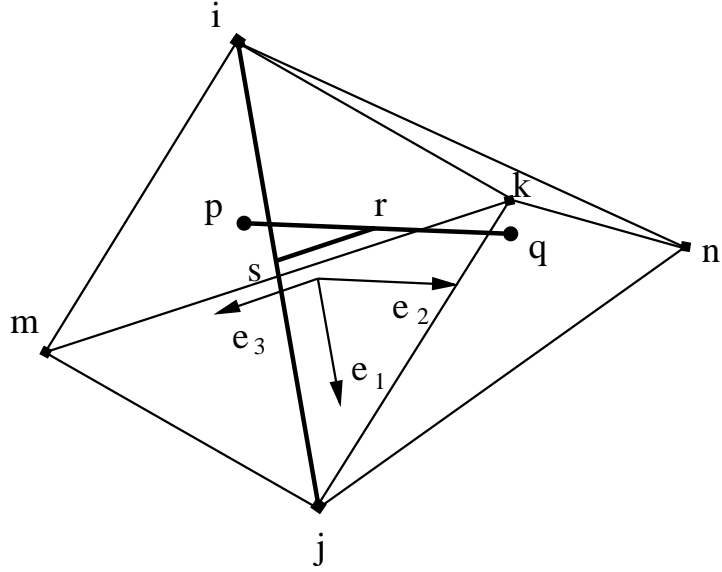} 
	\caption{Tree-dimensional case}
    \label{f2}
  \end{center}
\end{figure}

Here $i,j,k,m$ --- nodes of the Dalaunay triangulation, $i,j,k$ è $i,j,m$ --- 
two adjacent triangles of this triangulation with common edge $(i,j)$,  
$p$ è $q$ --- nodes of the Voronoi diagram.  For nodes of the Voronoi 
diagram there are used notations $\mathbf{x}^V_{p} = \mathbf{x}^D_{ijk}$ 
and $\mathbf{x}^V_{q} = \mathbf{x}^D_{ijm}$. Highlighted in Fig.\ref{f1}
edges of the Delaunay triangulation and Voronoi diagram are orthogonal to
one another. Introduce at the point of their intersection 
(at point $\mathbf{x} = \mathbf{x}^D_{ij}$) local coordinate system
$(s_1,s_2)$ with direction unit vectors $\mathbf{e}_1, \mathbf{e}_2)$.

Divergence of a vector field at the discrete level (see (\ref{11})) 
is defined using the tangential component. To approximate the diffusion
transport, it is necessary to evaluate the diffusion tensor (\ref{18}).
Assume for the introduced local coordinate system
\[
  (\mathbf{D})_{ij}^D = \left (
\begin{array}{cc}
 (d_{11}^D)_{ij}  &   (d_{12}^D)_{ij}  \\
 (d_{21}^D)_{ij}   &  (d_{22}^D)_{ij}   \\
\end{array}
\right ) .
\]
Coefficients $(d_{\alpha \beta }^D)_{ij}, \ \alpha, \beta =1,2$ are
calculated via specified smooth enough coefficients
$d_{\alpha \beta } (\mathbf{x}), \ \alpha, \beta =1,2$
of equation (\ref{1}). If local coordinate system $(s_1,s_2)$ coincides with 
$(x_1,x_2)$ we have, for instance, 
\[
  (d_{\alpha \beta }^D)_{ij}= 
  d_{\alpha \beta }(\mathbf{x}_{ij}^D),
  \quad  \alpha, \beta =1,2 .
\]
In the general case the recalculation rule for a tensor of second order
should be taken into account when going to a new orthogonal coordinate system.

Using the above notations in local coordinate system $(s_1,s_2)$ 
the tangential component of the flow can be written as
\begin{equation}\label{19}
  (\mathbf{D}\grad_D y)_{ij}^D  =
  (d_{11}^D)_{ij} (\grad_D y)_1 + 
  (d_{12}^D)_{ij} (\grad_D y)_2 .   
\end{equation}
The tangential component of the gradient seems like this
\begin{equation}\label{20}
(\grad_D y)_1 = 
  (\mathrm{grad}_D y)^D_{ij} =  \frac{y_j^D - y_i^D}{l_{ij}^D} .
\end{equation}
For the normal component we have
\begin{equation}\label{21}
(\grad_D y)_2 = 
  \frac{y_{q}^V - y_{p}^V}
  {l_{pq}^V},
\end{equation}
where $l_{pq}^V = |\mathbf{x}_{ijm}^D - \mathbf{x}_{ijk}^D |$ ---
the edge length of the Voronoi polygon (polyhedron).

Evaluation of the normal component of the gradient via (\ref{21})  is based on
interpolaton of values of a scalar function defined at the Delaunay
triangulation nodes to the Voronoi polygon nodes.
Let $l_{pij} = |\mathbf{x}_{p}^V - \mathbf{x}_{ij}^D | $ be the distance from 
node $p$ of the Voronoi partitioning up to edge $(i,j)$ 
of the Delaunay triangulation so that
$l_{pq}^V = d_{pij} + d_{qij}$. The linear interpolation with accuracy 
$\mathcal{O}(h^2)$ over triangle $D_p$ results in
\begin{equation}\label{22}
  y_{p}^V = \frac{1}{\mathrm{meas} (D_p)}
  \left ( y_k^D l_{ij}^D d_{pij} +
  y_i^D l_{jk}^D d_{pjk} +
  y_j^D l_{ik}^D d_{pik} \right ) .
\end{equation}
Taking into account (\ref{20}) - (\ref{22}),  
the next representation for the truncation error of the flow 
follows from (\ref{19})
\begin{equation}\label{23}
  (\mathbf{D}\grad_D u)_{ij}^D  =
  (\mathbf{D}\grad \,u)(\mathbf{x}_{ij}^D) +  \mathcal{O}(h) .
\end{equation}

Approximation (\ref{1}) leads to the grid equation
\[
  - \div_D (\mathbf{D}^D\grad_D y) + 
  \mathcal{C}_2^D (\mathbf{a}^D) y +
  \mathcal{C}_1^D (\mathbf{b}^D) y +
  c^D y = f^D,
\]
\begin{equation}\label{24}
  \quad \mathbf{x} = \mathbf{x}_{i}^{D},
  \quad i = 1,2,\ldots,M_D ,
\end{equation}
where, for instance, $c^D (\mathbf{x}_{i}^{D}) = c^D (\mathbf{x}_{i}^{D})$.
Grid equation (\ref{24}) approximates equation (\ref{1}) with error
\begin{equation}\label{25}
  \psi (\mathbf{x}) = g(\mathbf{x}) +
  (\mathrm{div}_D \mathbf{q}) ({\mathbf{x}}),
  \quad g = \mathcal{O}(h), 
  \quad \mathbf{q} = \mathcal{O}(h).   
\end{equation}

In the 3D case approximations of the diffusion transport operator
are constructed in a similar way. Let us consider (see Fig.\ref{f2})
two adjacent tetrahedrons of the Delaunay triangulation  --- 
$i,j,k,m$ and  $i,j,k,n$. For the corresponding nodes of the Voronoi diagram
(vortices of the Voronoi polygon) we again use notations
$\mathbf{x}^V_{p}$ and $\mathbf{x}^V_{q}$. The line connecting 
these vortices is perpendicular to triangle $i,j,k$ and intersects it at point
$r$ (Fig.\ref{f2}). Point $s$ is the midpoint of edge $(i,j)$. 
Stright line $(r,s)$ is perpendicular to both edge $(i,j)$
of the Delaunay triangulation and edge $(p,q)$ of the Voronoi
diagram. On the basis of these three lines $(i,j)$, $(p,q)$ and $(r,s)$ 
highlighted in Fig.\ref{f2} we construct local coordinate system
$(s_1,s_2,s_3)$ with direction unit vectors
$\mathbf{e}_1, \mathbf{e}_2, \mathbf{e}_3)$.

Similarly to the 2D case we can evaluate now the diffusion tensor 
in the local coordinate system
\[
  (\mathbf{D})_{ij}^D = \left (
\begin{array}{ccc}
 (d_{11}^D)_{ij}  &   (d_{12}^D)_{ij}  &   (d_{13}^D)_{ij} \\
 (d_{21}^D)_{ij}  &   (d_{22}^D)_{ij}  &   (d_{23}^D)_{ij} \\
 (d_{31}^D)_{ij}  &   (d_{32}^D)_{ij}  &   (d_{33}^D)_{ij} \\
\end{array}
\right ) .
\]
For the smooth enough coefficients
$d_{\alpha \beta } (\mathbf{x}), \ \alpha, \beta =1,2,3$ of equation (\ref{1})
and in case of coincidence of local coordinate system
$(s_1,s_2,s_3)$ with $(x_1,x_2,x_3)$ it is possible to take
\[
  (d_{\alpha \beta }^D)_{ij}= 
  d_{\alpha \beta }(\mathbf{x}_{ij}^D),
  \quad  \alpha, \beta =1,2,3 .
\]

In local coordinate system $(s_1,s_2,s_3)$ the flow component
along edge $(i,j)$ of the Delaunay triangulation is equal to
\[
  (\mathbf{D}\grad_D y)_{ij}^D  =
\]
\begin{equation}\label{26}
  (d_{11}^D)_{ij} (\grad_D y)_1 + 
  (d_{12}^D)_{ij} (\grad_D y)_2 + 
  (d_{13}^D)_{ij} (\grad_D y)_3 .   
\end{equation}
The tangential component of gradient $(\grad_D y)_1$ is calculated
in accordance with (\ref{20}) whereas for component $(\grad_D y)_2$
(along the edge of the Voronoi diagram) approximation (\ref{21})
is used. For the third term in (\ref{26}) we employ
\begin{equation}\label{27}
(\grad_D y)_3 = 
  \frac{\bar{y}_{r} - \bar{y}_{s}}
  {l_{rs}},
\end{equation}
where $l_{rs}$ --- distance between points $r$ and $s$.
Values of scalar function $\bar{y}_{r}, \bar{y}_{s}$ in (\ref{27})
are approximated with the second order via (\ref{22})-type formulas
using values defined at nodes of the Delaunay triangulation.

So, representation (\ref{23}) for the truncation error of the flow
is valid in the 3D case, too. Hence we aso obtain grid equation (\ref{24})
with representation (\ref{25}) for the truncation error.
On the basis of this divergent representation for the error
we can prove convergence of the discrete solution to the exact one
\cite{vago} in the grid analog of space $\mathcal{H}(\grad)$
with the first order with respect to $h$.

\bibliographystyle{elsarticle-num}
\bibliography{vago}

\begin{thebibliography}{10}
\expandafter\ifx\csname url\endcsname\relax
  \def\url#1{\texttt{#1}}\fi
\expandafter\ifx\csname urlprefix\endcsname\relax\def\urlprefix{URL }\fi
\expandafter\ifx\csname href\endcsname\relax
  \def\href#1#2{#2} \def\path#1{#1}\fi

\bibitem{Samarskii_theory}
A.~Samarskii, {The theory of difference schemes.}, Pure and Applied
  Mathematics, Marcel Dekker. 240. New York, NY: Marcel Dekker. 786 p., 2001.

\bibitem{balansametod}
A.~A. Samarskii, Parabolic equations and difference methods for their solution,
  in: Proc. of All-Union Conference on Differential Equations, 1958, Publishers
  of Armenian Ac.Sci., 1960, pp. 148--160.

\bibitem{FVM}
R.~Eymard, T.~Gallou\"{e}t, R.~Herbin, Finite volume methods, in: Handbook of
  Numerical Analysis, Vol.~7, North Holland, Amsterdam, 2000, pp. 713--1020.

\bibitem{Finite_volume_methods}
R.~J. Leveque, {Finite volume methods for hyperbolic problems.}, Cambridge
  Texts in Applied Mathematics. Cambridge: Cambridge University Press. xix, 558
  p., 2002.

\bibitem{Generalized_difference_methods}
R.~Li, Z.~Chen, W.~Wu, {Generalized difference methods for differential
  equations: Numerical analysis of finite volume methods.}, Pure and Applied
  Mathematics, Marcel Dekker. 226. New York, NY: Marcel Dekker. xv, 442 p.,
  2000.

\bibitem{George_Louis_Borouchaki}
P.-L. George, H.~Borouchaki, {Delaunay triangulation and meshing. Application
  to finite elements. Transl. from the French.}, {Paris: \'Edition Herm\`es.
  vii, 413 p.}, 1998.

\bibitem{DirichleFavorskii}
A.~Solov'ev, E.~Solov'eva, V.~Tishkin, A.~Favorskij, M.~Shashkov,
  {Approximation of finite-difference operators on a mesh of Dirichlet cells.},
  Differ. Equations 22 (1986) 863--872.

\bibitem{Nicolaides_Wang}
R.~Nicolaides, D.-Q. Wang, {Convergence analysis of a covolume scheme for
  Maxwell's equations in three dimensions.}, Math. Comput. 67~(223) (1998)
  947--963.

\bibitem{Nicolaides_Wu}
R.~A. Nicolaides, X.~Wu, {Covolume solutions of three-dimensional div-curl
  equations.}, SIAM J. Numer. Anal. 34~(6) (1997) 2195--2203.

\bibitem{vago}
P.~Vabishchevich, {Finite-difference approximation of mathematical physics
  problems on irregular grids.}, Comput. Methods Appl. Math. 5~(3) (2005)
  294--330.

\bibitem{Mishev}
I.~D. Mishev, Finite volume methods on voronoi meshes, Numerical Methods for
  Partial Differential Equations 14~(2) (1998) 193--212.

\bibitem{Vabishchevich_Samarskij2000}
P.~Vabishchevich, A.~Samarskij, {Finite difference schemes for
  convection-diffusion problems on irregular meshes.}, Comput. Math. Math.
  Phys. 40~(5) (2000) 692--704.

\bibitem{Lazarov_Vassilevski}
R.~D. Lazarov, P.~S. Vassilevski, {Numerical methods for convection-diffusion
  problems on general grids.}, {Bojanov, B. D. (ed.), Approximation theory. A
  volume dedicated to Blagovest Sendov. Sofia: DARBA. 258-283} (2002).

\bibitem{samvabkon}
A.~A. Samarskii, P.~N. Vabishchevich, Numerical Methods for Solving
  Convection-Diffusion Problems, Editorial URSS, Moscow, 2003, in Russian.

\bibitem{Matus_Rybak}
P.~Matus, I.~Rybak, {Difference schemes for elliptic equations with mixed
  derivatives.}, Comput. Methods Appl. Math. 4~(4) (2004) 494--505.

\bibitem{Putti_Cordes}
M.~Putti, C.~Cordes, {Finite element approximation of the diffusion operator on
  tetrahedra.}, SIAM J. Sci. Comput. 19~(4) (1998) 1154--1168.

\bibitem{Droniou_Eymard}
J.~Droniou, R.~Eymard, {A mixed finite volume scheme for anisotropic diffusion
  problems on any grid.}, Numer. Math. 105~(1) (2006) 35--71.

\bibitem{Eymard_Gallouet_Herbin}
R.~Eymard, T.~Gallou{\"e}t, R.~Herbin, {A cell-centred finite-volume
  approximation for anisotropic diffusion operators on unstructured meshes in
  any space dimension.}, IMA J. Numer. Anal. 26~(2) (2006) 326--353.

\bibitem{Hu_Nicolaides}
X.~Hu, R.~Nicolaides, {Covolume techniques for anisotropic media.}, Numer.
  Math. 61~(2) (1992) 215--234.

\bibitem{Gilbarg}
D.~Gilbarg, N.~S. Trudinger, Elliptic partial differential equations of second
  order, Grundlehren der Mathematischen Wissenschaften, 224. Berlin etc.:
  Springer-Verlag. XIII, 513 p., 1983.

\bibitem{PGCiarlet_1978a}
P.~G. Ciarlet, The Finite Element Method for Elliptic Problems, North--Holland,
  Amsterdam, New York, 1978.

\end{thebibliography}

\end{document}